\documentclass[a4paper]{article}
\usepackage{xparse}
\usepackage{amsmath}
\usepackage{amssymb}
\usepackage[english]{babel}
\usepackage[utf8]{inputenc}
\usepackage[]{nicefrac}
\usepackage{amsmath, amssymb, amsthm}
\usepackage{graphicx}
\usepackage[colorinlistoftodos]{todonotes}
\usepackage{tikz}
\usepackage[affil-it]{authblk}
\usepackage[]{ccicons}
\usepackage[]{xcolor}
\usepackage{mathtools}
\usepackage{enumitem}
\usepackage{appendix}
\usepackage{leftidx}
\usepackage[]{hyperref}
\hypersetup{
    colorlinks = true,
    urlcolor =gray}

\newcommand{\dVblank}{\, {\rm d}V}

\newcommand{\norm}[1]{\left\lVert#1\right\rVert}

\NewDocumentCommand \dV{ o }{%
    \IfNoValueTF{#1}{\dVblank}%
    {
        \dV_{#1}
    }
}

\NewDocumentCommand \mc{ m }{%
    \mathcal{#1}%
}

\NewDocumentCommand \eref{ m }{%
    (\ref{eqn:#1})%
}

\NewDocumentCommand \an{ m }{%
    \langle {#1} \rangle%
}

\NewDocumentCommand \mbf{ m }{%
    \mathbf{#1}
}

\NewDocumentCommand \stateOne{ m }{
    \underline{\mbf{#1}}%
}

\NewDocumentCommand \stateTwo{ m o }{%
    \IfNoValueTF{#2}{\stateOne{#1}}
    {
    \stateOne{#1}\an{#2}%
    }
}

\NewDocumentCommand \s{ m o o }{%
    \IfNoValueTF{#3}{\stateTwo{#1}[#2]}%
    {
    \underline{\mbf{#1}}({#3})\an{#2}%
    }
}

\NewDocumentCommand \stateOneI{ m m }{
    \underline{#1}_{#2}%
}

\NewDocumentCommand \stateTwoI{ m m o }{%
    \IfNoValueTF{#3}{\stateOneI{#1}{#2}}
    {
    \stateOneI{#1}{#2}\an{#3}%
    }
}

\NewDocumentCommand \si{ m m o o }{%
    \IfNoValueTF{#4}{\stateTwoI{#1}{#2}[#3]}%
    {
        \stateOneI{#1}{#2}(\mbf{#4})\an{#3}
    }
}

\definecolor{officegreen}{rgb}{0.0, 0.5, 0.0}

\title{Machine-learning of nonlocal kernels\\ for anomalous subsurface transport\\ from breakthrough curves}
\date{\today}

\author{Xiao Xu\thanks{The Oden Institute for Computational Engineering and Sciences, The University of Texas at Austin, TX, xiaoxu42@utexas.edu}
\and Marta D'Elia \thanks{Computational Science and Analysis, Sandia National Laboratories, CA, mdelia@sandia.gov}
\and Christian Glusa \thanks{Computer Science Research Institute, Sandia National Laboratories, NM, caglusa@sandia.gov}
\and John T. Foster\thanks{The Oden Institute for Computational Engineering and Sciences, The University of Texas at Austin, TX, john.foster@utexas.edu}
}

\begin{document}
\maketitle
\begin{abstract}
    Anomalous behavior is ubiquitous in subsurface solute transport due to the presence of high degrees of heterogeneity at different scales in the media. Although fractional models have been extensively used to describe the anomalous transport in various subsurface applications, their application is hindered by computational challenges. Simpler nonlocal models characterized by integrable kernels and finite interaction length represent a computationally feasible alternative to fractional models; yet, the informed choice of their kernel functions still remains an open problem. We propose a general data-driven framework for the discovery of optimal kernels on the basis of very small and sparse data sets in the context of anomalous subsurface transport. Using spatially sparse breakthrough curves recovered from fine-scale particle-density simulations, we learn the best coarse-scale nonlocal model using a nonlocal operator regression technique. Predictions of the breakthrough curves obtained using the optimal nonlocal model show good agreement with fine-scale simulation results even at locations and time intervals different from the ones used to train the kernel, confirming the excellent generalization properties of the proposed algorithm. A comparison with trained classical models and with black-box deep neural networks confirms the superiority of the predictive capability of the proposed model. 
\end{abstract}
\section{Introduction}

Making accurate large-scale predictions of solute transport in the subsurface is critically important for the efficient management of water resources \cite{Sun2020,suzuki2021fractional}
as well as petroleum production, particularly in enhanced oil-recovery (EOR) applications \cite{dehghani1984, srikanta1991, stalgorova2012, tong2017}. Subsurface transport is a highly complex phenomenon as it takes place in environments that contain heterogeneities at all scales, requiring the use of fine-scale models at the smallest scales. However, it is well-known that direct numerical simulations of subsurface transport at the fine scales can be prohibitively expensive, despite the recent advances in computational power and high-performance computing. This creates the need for upscaled (or coarse-grained) models that act at large scales and are able to capture the effects of the fine-scale behavior. An open challenge is the fact that, due to the presence of high degrees of heterogeneties, the models that are accurate at the fine scales fail to be predictive at large scales. In other words, at coarser scales, quantities of interest (such as the particle concentration) follow laws that are substantially different from the ones considered accurate at the fine scales \cite{neuman2009perspective}.  Another ongoing challenge is the fact that is it nearly impossible to accurately measure relevant medium properties, so that, even when a model is available, its parameters remain uncertain. Despite these difficulties, several works have addressed the problem of simulating particle transport at large scales, while taking into account the effects of the fine-scale behavior. Relevant to this work, we mention the use of nonlocal models: it has been shown (see, e.g., \cite{levy2003measurement}) that heterogeneities at small scales result in non-Fickian (anomalous) behavior when upscaled to larger scales; specifically, they yield a nonlinear mean squared displacement (MSD), as opposed to the classical, linear MSD. These anomalous effects often take the form of early arrival times and late-time long tails in breakthrough curves (BTCs), i.e. concentration profiles at a specific location as a function of the time.
In recent years, substantial effort has been devoted to the development of fractional advection-dispersion equations for describing anomalous transport in heterogeneous environments \cite{d2021fractional,d2021analysis,Kelly2019chapter,Schumer2001,suzuki2021fractional}, as these models intrinsically yield the desired nonlinear MSD and embed all length scales in their definition. At the same time, fractional models come with several computational challenges, mostly due to their high computational cost and to their nontrivial implementation. Recent developments in nonlocal theory have shown that fractional operators are special cases of more general nonlocal operators \cite{d2021connections,DElia2021Unified,d2013fractional}; specifically, for certain choices of kernel functions, fractional operators are the limits of nonlocal operators as the extent of the nonlocal interactions goes to infinity. These facts motivate the model that we propose in this work:
{\it we conjecture that a spatial nonlocal operator featuring long-range interactions of finite length provides a reliable description of the coarse-grained evolution of the particle density in heterogeneous subsurface environment.}

The simplest and most general form of nonlocal operator $\mc L_\omega$ applied to a scalar function $u$ is given by \cite{DElia2021Unified,Du2012}
\begin{align}
    \mc L_\omega u(x) = \int_{\mc H} \omega(\xi) (u(x+\xi) - u(x)) d \xi, \notag 
\end{align}
where $\omega$ is the nonlocal kernel function, and $\mc H$ is a neighborhood of $x$ of size $\delta>0$. The latter, often referred to as the \emph{horizon}, determines the extent of the nonlocal interactions. While the kernel function $\omega$ is a fundamental entity that determines the intrinsic behavior of the quantity of interest $u$ and its regularity properties ~\cite{weckner2005effect,seleson2011role}, {\it there is still no general technique to systematically determine such kernel a priori.}

Recently, a few attempts have been made to derive the nonlocal kernel function for nonlocal models in the context of solid mechanics applications. In \cite{silling2014origin}, it is demonstrated that nonlocality can arise from the micro-scale heterogeneity as the effect of an implicit or explicit homogenization procedure. In \cite{xu2020deriving}, a theoretical method is proposed to determine the kernel function based on the micro-structure of a periodic, heterogeneous, one-dimensional bar. In \cite{you2020data,xu2021machine,You2021AAAI}
the authors use machine learning to identify optimal kernel function for stress wave propagation from synthetic high-fidelity data. More recently, in \cite{You2021MD}, the same machine learning techniques have been used to learn nonlocal kernels with the purpose of reproducing molecular dynamics coarse-grained behavior.

Inspired by these works, we propose to use similar techniques to derive data-driven, coarse-grained, nonlocal models that describe anomalous transport in heterogeneous subsurface environments. With the purpose of reproducing the nonlinear MSD, we introduce the time-dependence of the kernel function and refer to the resulting model as \emph{nonlocal dynamic kernel}. Motivated by the fact that BTCs are usually the available field data, we infer the optimal kernel function solely on the basis of spatially sparse BTCs over a limited time interval. In this work, for simplicity and accurate validation, we generate the BTCs via high-fidelity fine-scale simulations and leave more realistic data sets to future works. The curves used for training are the result of coarse-graining of the fine-scale particle density. We summarize our main contributions below.

\begin{itemize}
    \item We design a general machine learning framework that uses solely BTC data in order to learn nonlocal operators that describe anomalous diffusion.
    \item We add a time-dependent factor to the nonlocal kernel that allows us to capture the nonlinear MSD of the anomalous diffusion.
    \item We show that the model learned from early-time BTC data can generalize well and make accurate predictions of the late-time long tails of the BTC.
    \item We also show the improved performance of our learned model compared to classical PDE models and surrogates based on artificial neural networks (NN).
\end{itemize}

\paragraph{Outline of the paper} In \S\ref{chap:local}, we introduce the two-dimensional periodic heterogeneous porous medium, the fine-scale, and the high-fidelity governing equations. In \S\ref{chap:nonlocal}, we introduce the upscaling (or coarse-graining) procedure, and describe the proposed nonlocal diffusion model for the coarse-grained particle density. \S\ref{chap:learning} outlines the learning algorithm and its discretization. In \S\ref{chap:numerical} we report several computational tests that illustrate the generalization properties of the proposed approach as well as comparisons with state-of-the-art models and deep learning approaches. \S\ref{chap:conclusion} summarizes our contributions and provides future research guidelines.


\section{A high-fidelity model for the particle density}
\label{chap:local}

In this section we introduce the high-fidelity model that will be used for the generation of the coarse-grained BTCs and for comparison of the predicted curves with out-of-training-range simulations. We first describe the equations governing the flow through heterogeneous media and then introduce the equations describing the evolution of the particle density. 

\subsection{Flow through a heterogeneous medium}
\begin{figure}
\centering
\begin{tikzpicture}[scale=0.65]
    \tikzstyle{pt}=[circle, fill=black, inner sep=0pt, minimum size=5pt]
    \draw (-8,0) -- (8,0);
    \draw (8,0.8) -- (-8,0.8);
    \draw (8,0) -- (8,0.8);
    \draw (-8,0) -- (-8,0.8);
    \foreach \x in {-2.4-4,-1.6-4,-0.8-4,0-4,0.8-4,1.6-4,2.4-4,3.2-4,0,0.8,1.6,2.4,3.2,4.0,4.8,5.6,6.4}{
        \node[rectangle, draw, minimum height=0.6, fill=black!50, yscale=1.2, rotate=45] at (\x,0.4) {};
    }

    
    
    \draw[->] (-3.4,0.4)  to [out=220,in=160, looseness=1] (-4.2,-2.8);

	\draw[-latex] (-8-0.5,-0.4-0.5) -- (-6-0.5,-0.4-0.5); 	
	\draw[-latex] (-8-0.5,-0.4-0.5) -- (-8-0.5,1.2-0.5); 	
    \draw (-8,-0.7-0.5) node() {$x$};
    \draw (-8.3,1.2-0.5) node() {$y$};
    
    \draw[-latex] (-3.2,2.4) -- (-3.2,0.8) node[near start, above , fill=white, inner sep=1pt] {$\textrm{particle injection}$} ;
    \draw[latex-latex] (-8,1.2) -- (-3.2,1.2) node[midway , above,  fill=white, inner sep=1pt] {$L_1$} ;
    
    \draw[latex-latex] (8,-0.5) -- (-8,-0.5) node[midway , fill=white, inner sep=1pt] {$L$} ;
    
   
   \filldraw[fill=white, draw=black] (-4,-2.5) rectangle (0, -6.5);
   
   \node[rectangle, draw, fill=black!50,xscale=1*6, yscale=1.2*6, rotate=45] at (-2,-4.5) {};
    
    \draw (-4,-2.5) -- (-4,-1.9);
    
    \draw (0,-2.5) -- (0,-1.9);
    
    \draw[latex-latex] (-4,-2.2) -- (0,-2.2) node[midway , above,  fill=white, inner sep=1pt] {$l_1$} ;
    
    \draw[latex-latex] (0.3,-2.5) -- (0.3,-6.5) node[midway , above,  fill=white, inner sep=1pt] {$l_2$} ;
    \draw (0,-2.5) -- (0.6,-2.5);
    \draw (0,-6.5) -- (0.6,-6.5);
    
    \filldraw[fill=black!50, draw=black] (2,-2-1) rectangle (2.8,-2.8-1);
    \draw (4,-2.4-1) node() {$ \kappa$};

    \filldraw[fill=white, draw=black] (2,-3-1) rectangle (2.8,-3.8-1);
    \draw (4,-3.4-1) node() {$ \kappa_0$};

    

\end{tikzpicture}
\caption{Sketch of a two-dimensional porous medium with periodic heterogeneity.}
\label{fig:2DIllustration}
\end{figure}
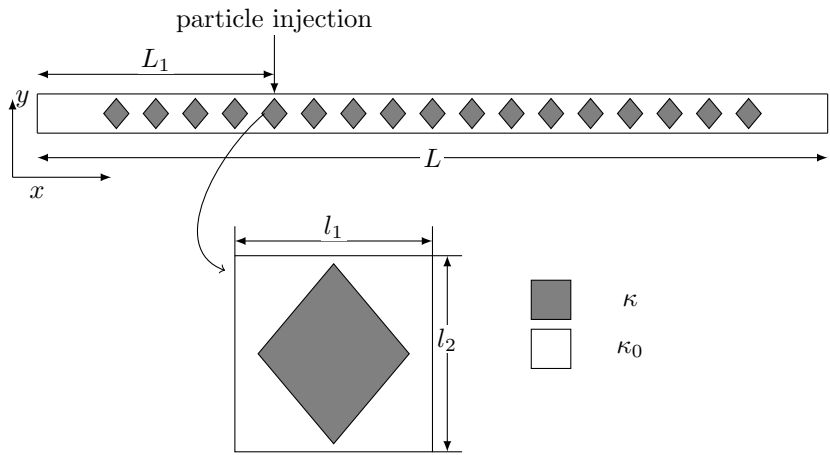

We consider the flow field through a heterogeneous porous medium, following the same setting as in \cite{tyukhova2016mechanisms}. We refer to the simplified configuration sketched in Figure~\ref{fig:2DIllustration}, where the porous medium is a two-dimensional thin layer composed of a homogeneous matrix of conductivity $\kappa_0$ and periodic diamond-shaped inclusions of conductivity $\kappa$. Let the thickness of the layer be $l_2$ and the length of the unit cell of an inclusion be $l_1$, and assume the length-scale of the unit cell is much smaller than the horizontal length-scale of the layer $L$,
\begin{align}
    l_1,l_2 \ll L , \notag 
\end{align}
Let $N$ denote the total number of unit cells in the medium, so that, for consistency,
\begin{align}
    L = N l_1. \notag
\end{align}
We assume a constant hydraulic head (i.e. liquid pressure) $h_0$ on the left boundary, a zero hydraulic head on the right boundary, and no flow across the top and bottom boundaries. Furthermore, we assume that there is no flow source inside the medium. Then, the flow field $\mbf v(x,y)$ and the hydraulic head $h(x,y)$ inside the two-dimensional domain $\Omega = [0,L]\times[0,l_2]$ satisfy the Darcy equations
\begin{equation}
    \begin{cases}
        \mbf v(x,y) = - \kappa(x,y) \nabla h(x,y) \\
        \nabla \cdot \mbf v(x,y) = 0 
    \end{cases} \quad (x,y) \in \Omega
    \label{eqn:Poieqn} 
\end{equation}
with boundary conditions
\begin{align}
  \left\{
  \begin{array}{rll}
    h(0,y) &= h_0, & y\in(0,l_{2}), \\
    h(L,y) &= 0, & y\in(0,l_{2}), \\
    v_y(x,0) &= 0, & x\in(0,L), \\
    v_y(x,l_2) &= 0, & x\in(0,L).
  \end{array}
  \right. 
  \label{eqn:PoieqnBC} 
\end{align}

%

\subsection{Particle transport in heterogeneous medium}
We consider the transport of particles injected in the medium induced by the advection field $\bf v$, introduced in the previous section, in the absence of fine-scale diffusion.
Let $c(x,y,t)$ denote the density of the particles at $(x,y)$ and time $t$. Then the evolution of the particle density $c(x,y,t)$ induced by flow field $\mbf v(x,y)$ satisfies the Liouville equation
\begin{align}
    \frac{\partial c(x,y,t)}{\partial t} = -\mbf v(x,y) \cdot \nabla c(x,y,t). \notag 
\end{align}
Equivalently, using a Lagrangian viewpoint, the evolution of the position of each particle can be described by the following equation:
\begin{align}
    \frac{d \mbf x(t)}{dt} = \mbf v(\mbf x(t)),
    \label{eqn:Lageqn}
\end{align}
where $\mbf x$ is the coordinate vector $\mbf x = (x,y)^T$ of the particle. We refer to \eqref{eqn:Lageqn} as the particle dynamic equation. Under the assumption that the particles are injected into the \(k\)-th unit cell of the medium, \(C_{k}:=[(k-1)l_{1},k l_{1}]\times[0,l_{2}]\), proportionally to the flux inside that unit cell at time $t=0$, we have the following initial conditions:
  \begin{equation*}
    c(x,y,0)=
  \begin{cases}
    |\mbf v(x,y)| & \text{ if } (x,y)\in C_{k} , \\
    0 & \text{otherwise.}
  \end{cases}
  \end{equation*}
Despite their prohibitive computational cost in realistic settings, discretizations of \eqref{eqn:Lageqn} for a large number of particles are often used in the literature to describe the transport of particles in heterogeneous subsurface environments, provided an advection field is available \cite{edery2014origins, tyukhova2016mechanisms}. As such, we will treat this model and its discretization as the high-fidelity reference to be used for data generation and comparison.


\section{The proposed nonlocal model and its discretization}
\label{chap:nonlocal}

In this section we describe the proposed nonlocal model for the evolution of the coarse-grained particle density that is used to predict the evolution of the BTCs. First, we introduce the coarse-graining technique and then propose a coarse-grained nonlocal model, characterized by a time-dependent kernel whose choice is motivated by the anomalous nature of the diffusive process. We also recall a well-established discretization technique for the nonlocal diffusion equation that will be used in our experiments. 

\subsection{A nonlocal upscaled equation for the particle density}
\label{sec:upscaling}

With the purpose of describing particle transport at scales that are much larger than the extent of the heterogeneities, following the claim made in the introduction, we propose to upscale the two-dimensional heterogeneous model for the particles' density to a one-dimensional model characterized by nonlocal operators. To this end, we introduce the coarse-scale particle density $\bar{c}(x, t)$, defined as the average of the fine-scale density $c(x,y,t)$ over the unit cells.
Specifically, for each $x \in [0,L)$, let $N_x$ be the unique integer such that $(N_x-1) l_1 \leq x < N_x l_1$. We define the coarse-scale particle density at $x$ as the average of the fine-scale particle density in $m$ unit cells, i.e., the average over $C_{N_k} \cup C_{N_k+1} \dots \cup C_{N_k+m-1}$:
\begin{align}
    \bar{c}(x, t) = \frac {\int_{0}^{l_2}\int_{(N_x-1) l_1}^{(N_x+m-1) l_1} c(\xi,\eta, t) d\xi d\eta }{ m l_1 l_2 }. 
    \label{eqn:upscaling}
\end{align}
The parameter $m$ is an integer that determines the smoothness of the coarse-scale density: the larger $m$, the smoother the density, as the average occurs over a larger set of cells.

{The rate of change of \(\bar{c}\) is determined by the flux across the interfaces \(x=(N_{x}-1)l_{1}\) and \(x=(N_{x}+m-1)l_{1}\).
Allowing for long-range jumps, we}
conjecture that the coarse-scale density $\bar{c}(x,t)$ satisfies a nonlocal parabolic equation of the following form:  
\begin{align}
    \frac{\partial \bar c(x,t)}{\partial t} = \int_{-\delta}^{\delta} \omega(\xi,t)( \bar c(x+\xi,t) - \bar c(x,t)  ) d\xi - \bar v \frac{\partial \bar c(x,t)}{\partial x},
    \label{eqn:PDeqn}
\end{align}
provided the nonlocal Dirichlet boundary condition
\begin{align}
    \bar c(x,t) = 0 \qquad x\in [-\delta,0) \cup (L,L+\delta] \notag 
\end{align}
is satisfied. Here, $\bar v$ is a constant scalar describing the ``effective'' coarse-scale advection and has the dimensions of a velocity field. In classical homogenization approaches, it is computed as:
\begin{align}
    \bar v = \frac{h_0}{N\bar\kappa_x}, \notag 
\end{align}
where $\bar \kappa_x$ is the effective conductivity of the heterogeneous unit cell in the $x$-direction (for its detailed derivation, see Appendix \ref{app:avgv}). Equation \eqref{eqn:PDeqn} is an advection-diffusion equation where the advection term is classical and the diffusion term is a nonlocal Laplacian with a time-dependent kernel $\omega$. We refer to the latter as nonlocal ``dynamic'' kernel and provide more details on its choice in Section \ref{sec:dynamic-kernel}. 

Due to the established relationship, in the case of fractional models, between the nonlocal kernel and the jump rate of the stochastic process associated with the fractional equation \cite{d2017nonlocal,suzuki2021fractional}, we assume that the nonlocal kernel is nonnegative. This assumption, enforced as a constraint in the learning procedure described in the following section, also guarantees the well-posedness of problem \eref{PDeqn} \cite{d2017nonlocal}.

\subsection{Separating diffusion from advection}
The classical, differential nature of the constant-advection term in \eref{PDeqn}, allows us to separate the contribution of the nonlocal diffusion and the classical advection terms by a simple change of variables.  Let $x_d = x - \bar v t$ and $\bar c_d (x_d,t) = \bar c(x,t)$; then, from \eref{PDeqn}, we see that $\bar c_d$ satisfies the following nonlocal diffusion equation without the advection term:
\begin{align}
    \frac{\partial \bar c_d(x_d,t)}{\partial t} = \int_{-\delta}^{\delta} \omega(\xi,t)( \bar c_d(x_d+\xi,t) - \bar c_d(x_d,t)  ) d\xi.
    \label{eqn:PDdiffeqn}
\end{align}
The main benefit of removing the classical advection term is to avoid numerical complications caused by the discretization of the advection term; this fact is essential for the efficiency of our learning algorithm. For the sake of notation, from now on, we drop the subscripts of $\bar c_d$, $x_d$.
After separating the nonlocal diffusion from the effective coarse-scale advection, the density $\bar c(x,t)$ governed by the diffusion-only nonlocal equation \eref{PDdiffeqn} {in an infinite domain} must satisfy the non-advective condition, i.e. the average displacement $\an x$ remains constant over time:
\begin{align}
    \frac{d \an x (t)}{dt} = 0,
    \label{eqn:nonadvective_cond}
\end{align}
with the average displacement $\an x$ given by
\begin{align}
    \an x (t)= \frac{\int^\infty_{-\infty} \bar c(x,t) x dx}{C},
    \label{eqn:averagedisp}
\end{align}
where $C:=\int^\infty_{-\infty} \bar c(x,t) dx$ is a constant representing the total number of particles inside the domain $[0,L]$. Substituting \eref{PDdiffeqn} into \eref{nonadvective_cond} and \eref{averagedisp}, we have
\begin{align}
    0=& \int^\infty_{-\infty} \frac{\partial \bar c(x,t)}{\partial t} x dx \notag \\
    =& \int^\infty_{-\infty}  \int_{-\delta}^{\delta}(\bar c(x+\xi,t)-\bar c(x,t))\omega(\xi,t)  x d\xi dx \notag \\
    =& \int_{-\delta}^{\delta}\int^\infty_{-\infty}  \bar c(x+\xi,t) \omega(\xi,t)  x dx d\xi -  \int_{-\delta}^{\delta}\int^\infty_{-\infty} \bar c(x,t)\omega(\xi,t)  x dx d\xi. \notag
\end{align}
By the change of variables $z=x-\xi$, we have  
\begin{align}
    =& \int_{-\delta}^{\delta}\int^\infty_{-\infty} \bar c(z,t)\omega(\xi,t) (z+\xi) dz d\xi -   \int_{-\delta}^{\delta}\int^\infty_{-\infty} \bar c(x,t)\omega(\xi,t)  x dx d\xi, \notag \\
    =& \int_{-\delta}^{\delta}\int^\infty_{-\infty} \bar c(z,t)\omega(\xi,t) (z+\xi) dz d\xi -   \int_{-\delta}^{\delta}\int^\infty_{-\infty} \bar c(x,t)\omega(\xi,t)  x dx d\xi \notag \\
    =& \int_{-\delta}^{\delta}\int^\infty_{-\infty} \bar c(z,t)\omega(\xi,t) \xi dz d\xi \notag \\ 
    =& C \int_{-\delta}^{\delta}\omega(\xi,t)\xi d\xi. \notag 
    \end{align}
Therefore, we have the non-advective constraint for the nonlocal diffusion kernel $\omega$:
\begin{align}
    \int^\delta_{-\delta} \omega(\xi,t) \xi d\xi = 0,
    \label{eqn:nonadvcon}
\end{align}
which ensures that the nonlocal diffusion kernel $\omega$ will not contribute to the average advection displacement $\an x$.

\subsection{The nonlocal dynamic kernel}\label{sec:dynamic-kernel}
We elaborate the importance of adding the time-dependency to the nonlocal kernel $\omega$ in relation to the mean square displacement of the particle density using similar approaches as described in the previous section. Consider the density $\bar c(x,t)$ satisfying the diffusion-only nonlocal equation \eref{PDdiffeqn} {in an infinite domain}. {The mean square displacement is defined by
\begin{align}
    \textrm{MSD} = \an{x - x_0}^2, \notag 
\end{align}
where $x_0$ is the initial position of the particles. The time derivative of the MSD is given by
\begin{align}
      \frac{d }{dt}\textrm{MSD} =& \frac{d \langle x^{2} - 2x_0 x + x_0^2 \rangle (t)}{dt} \notag \\
      =& \frac{d \langle x^{2}\rangle (t)}{dt} \notag \\
      =& \frac{1}{C}\frac{\partial(\int^\infty_{-\infty} \bar c(x,t)x^2 dx)}{\partial t},
    \label{eqn:MSD}
\end{align}
}
{where we used \eref{nonadvective_cond} to eliminate the drift term.} Substituting \eref{PDdiffeqn} into \eref{MSD}, we have 
\begin{align}
    C\frac{d }{dt}\textrm{MSD}=& \int^\infty_{-\infty} \frac{\partial \bar c(x,t)}{\partial t} x^2 dx \notag \\
    =& \int_{-\delta}^{\delta}\int^\infty_{-\infty}  \bar c(x+\xi,t) \omega(\xi,t)  x^2 dx d\xi -  \int_{-\delta}^{\delta}\int^\infty_{-\infty} \bar c(x,t)\omega(\xi,t)  x^2 dx d\xi. \notag \\
    =& \int_{-\delta}^{\delta}\int^\infty_{-\infty} \bar c(z,t)\omega(\xi,t) (z+\xi)^2 dz d\xi -   \int_{-\delta}^{\delta}\int^\infty_{-\infty} \bar c(x,t)\omega(\xi,t)  x^2 dx d\xi, \notag \\
    =& \int_{-\delta}^{\delta}\int^\infty_{-\infty} \bar c(z,t)\omega(\xi,t) (z+\xi)^2 dz d\xi -   \int_{-\delta}^{\delta}\int^\infty_{-\infty} \bar c(x,t)\omega(\xi,t)  x^2 dx d\xi \notag \\
    =& \int_{-\delta}^{\delta}\int^\infty_{-\infty} \bar c(z,t)\omega(\xi,t)(2z\xi + \xi^2) dz d\xi. \notag 
\end{align}
By using the non-advective constraint \eref{nonadvcon} on the first term we then conclude that
\begin{equation}\label{eqn:MSEproof}
C\frac{d }{dt}\textrm{MSD} = \int_{-\delta}^{\delta}\int^\infty_{-\infty} \bar c(z,t)\omega(\xi,t)\xi^2dz d\xi 
=C \int_{-\delta}^{\delta}\omega(\xi,t)\xi^2 d\xi. 
\end{equation}
The calculations above imply that when the kernel $\omega$ is independent of time, the mean square displacement of the density can only be linear function of time since its time-derivative is a constant. This would greatly limit the applicability of the proposed nonlocal model, as it would not be able to capture anomalous effects, i.e. a nonlinear behavior of the mean squared displacement. The latter is indeed achievable with the proposed, time-dependent kernel. With the only purpose of reducing the complexity of the learning problem, we simplify the expression of the kernel by separating the space and time dependencies, i.e.
\begin{equation}\label{eqn:omegasep}
    \omega(\xi,t) = \phi(\xi) \theta(t).
\end{equation}
where $\phi$ is the spatial kernel function and $\theta$ is the temporal kernel function. Substituting \eqref{eqn:omegasep} into \eqref{eqn:MSEproof} gives
\begin{align}
    \frac{d }{dt}\textrm{MSD} =  \int_{-\delta}^{\delta}\phi(\xi)\xi^2 d\xi \theta(t). \notag 
\end{align}
To further simplify the expression of the kernel and to mimic the power-law behavior of the mean squared displacement observed in the presence of anomalous diffusion, we select the following model for $\theta$:
\begin{align}
    \theta(t) = t^p, \notag 
\end{align}
where $p$ is a modeling parameter. With this choice, we conclude that the mean squared displacement is $\textrm{MSD} \sim t^{p+1}$.

\subsection{{Discretization}} 
\label{sec:space-discr}
Owing to definition \eref{upscaling}, the coarse-scale density $\bar{c}(x)$ is constant inside each unit cell. This fact motivates the use of a strong-form, piece-wise constant spatial discretization of the nonlocal diffusion equation where each unit cell corresponds to a degree of freedom. In this setting, we denote by $\bar c_i(t)$ the value of the coarse-scale density over the unit cell $C_i$, for $i=1,\ldots,N$. For the sake of simplicity, we assume that the horizon is a multiple of the length of unit cell, i.e.
\begin{align}
    \delta = N_\delta l_1, \notag 
\end{align}
where $N_\delta>0$ is an integer. Then, by using a Riemann sum over the neighboring unit cells to numerically approximate the integral, the discretized nonlocal diffusion equation equation reads \cite{silling2005meshfree}
\begin{align}
    \frac{\partial \bar c_i(t)}{\partial t} = \sum_{j=-N_\delta}^{j=N_\delta} \omega_j(t) (\bar c_{i+j}(t) - \bar c_i(t)  ) \quad 1\leq i\leq N
    \label{eqn:dis-PDdiffeqn}
\end{align}
where $\omega_j(t)$ is the discretized kernel function satisfying, for $-N_\delta\leq j\leq N_\delta$
\begin{align}
    \omega_j(t) =& \int_{(j-\frac{1}{2})l}^{(j+\frac{1}{2})l} \omega(\xi,t) d\xi \notag \\
    =& \theta(t) \int_{(j-\frac{1}{2})l}^{(j+\frac{1}{2})l} \phi(\xi) d\xi \notag \\
    =& \phi_j t^p.  \notag 
\end{align}

The temporal discretization of \eref{dis-PDdiffeqn} is done uniformly using first-order implicit backward differences.





\section{Kernel Learning}
\label{chap:learning}

In this section we describe the proposed learning algorithm for the kernel function $\omega(\xi,t)$. We assume that the only available high-fidelity data are sparse BTCs of the coarse-grained particle concentration through the heterogeneous medium described in Section \ref{chap:local}, i.e.
\[ \{ \left(f^*_i(t_j) , x_i \right) \vert \; i=1,2,\dots, n_x \quad j=1,2,\dots, n_t \}, \]
where the BTC at position $x_i$ and time $t_j$ is given by
\begin{equation}\label{eq:BTC}
        f^*_i(t_j) = \bar c_{H\!F}(x_i,t_j). 
    \end{equation}
Here, the subscript $H\!F$ stands for ``high fidelity'' and is added to distinguish the data from the learned density. While the index $j$ spans the whole time-discretization domain, the index $i$ spans a very small subset of the space-discretization domain, so that $n_x\ll N$. In particular, in our experiments, we will set $n_x\leq 3$. We also assume that $\{x_i\}_i$ are relatively far from the left side of the domain.
In general, the BTCs can be either obtained from field observations or generated by numerically solving the local equations in Section \ref{chap:local} with high accuracy. In our experiments, we will follow the latter approach. Below, we summarize the algorithmic workflow, from data generation, to kernel learning.

\noindent \underline{Data generation} We compute the fine-scale particle density $c(x,y,t)$ and upscale it to one-dimensional coarse-scale particle density $\bar c_{H\! F}(x,t)$ using \eref{upscaling}. The BTCs are then computed using \eqref{eq:BTC} at locations $\{x_i \vert \; i=1,2,\ldots n_x\}$.
    
\noindent \underline{Learning algorithm} We seek the optimal $\omega>0$ that minimizes the following loss function
    \begin{align}
        loss(\omega) = M\!S\!E(\omega) + \beta M(\omega), \notag 
    \end{align}
    where $M\!S\!E$ is the mean square error of the predicted BTCs
    \begin{align}
        M\!S\!E(\omega) = \sum_{i=1}^{n_x} \norm{f_i(\omega,\cdot) - f^*_i}^2_2, \notag 
    \end{align}
    where $\norm{\cdot}_2$ denotes the $L^2$ norm in the time domain, and $f_i(\omega,\cdot)$ is the BTC associated with the nonlocal, coarse-grained solution of \eref{PDeqn} for some kernel $\omega$ at a given location $x_i$.
    $M$ is the residual of the non-advective constraint $\eref{nonadvcon}$, i.e.
    \begin{equation}
        M(\omega) = \left(\int_{-\delta}^{\delta} \omega(\xi) \xi {\rm d}y \right)^2, \notag 
    \end{equation}
    and $\beta$ is a penalization parameter. 
    
\noindent\underline{Numerical solution}
For the discretization of the nonlocal equation we use the scheme described in Section \ref{sec:space-discr}. The set of unknowns that uniquely characterizes the kernel is given by the power $p$ and the vector $\nu_\phi \in \mathbb{R}^{2 N_\delta+1}$ defined as
\begin{align}
    \nu_\phi = (\phi_{-N_\delta},\phi_{1-N_\delta},\dots,\phi_{N_\delta-1},\phi_{N_\delta}). \notag 
\end{align}
We discretize each term in the loss function and, with an abuse of notation, we use the same symbols to indicate the discretized quantities. The discrete $M\!S\!E$ is given by
\begin{equation}
    M\!S\!E(\nu_\phi,p)= \sum_{i=1}^{n_x} \sum_{j=1}^{n_t} \left(f_i(\nu_\phi,p,t_j) - f^*_i(t_j)\right)^2. \notag 
\end{equation}
On the other hand, the discrete non-advective constraint is given by 
\begin{align}
    M(\nu_\phi) = \left( \sum^{j=N_\delta}_{j=-N_\delta} j \phi_j \right)^2. \notag 
\end{align}
The discretized loss function 
\begin{equation}\label{eq:discrete-loss}
   loss(\nu_\phi,p) = M\!S\!E(\nu_\phi,p) + \beta M(\nu_\phi)   
\end{equation}
is minimized with respect to $(\nu_\phi,p)$ using the L-BFGS algorithm with automatic differentiation {implemented in ForwardDiff.jl \cite{revels2016forward}}. The positivity constraint on the discrete $\omega$ is enforced during optimization by using a softplus function.

\section{Numerical Experiments}
\label{chap:numerical}

In this section we illustrate the efficacy of the proposed nonlocal model and the generalization properties of our learning strategy. We do this by comparing the predicted BTCs with high-fidelity simulations obtained with the model described in Section \ref{chap:local} in correspondence of locations and time intervals different from the ones used for training. Furthermore, to evaluate the performance of the proposed diffusion model, we compare it with three models that could be considered valid alternatives to \eqref{eqn:PDdiffeqn}. The first two are coarse-grained local models. Specifically, the first model uses the so-called ``fractal derivative'' \cite{CHEN20101754}
\begin{align}
    \frac{\partial \bar c(x,t)}{\partial t} = \frac{\bar D}{t^q} \frac{\partial^2 \bar c(x,t)}{\partial x^2},
    \label{eqn:fractaladvdiffeqn}
\end{align}
where, following the same arguments used to define the dynamic nonlocal kernel, we consider an effective diffusivity $\bar D$ scaled by the $q$-th power of the time. The second model is the classical diffusion model
\begin{align}
    \frac{\partial \bar c(x,t)}{\partial t} = \bar D_0 \frac{\partial^2 \bar c(x,t)}{\partial x^2} 
    \label{eqn:advdiffeqn}
\end{align}
where $\bar D_0$ is the effective diffusivity. To allow for fair comparisons, the parameters $D$, $q$, and $D_0$ are identified using the same learning algorithm described in Section \ref{chap:learning}. The third model is not a PDE and does not rely on physical arguments: we consider the approximation of the BTC via neural networks (NNs) {implemented in Pytorch \cite{NEURIPS2019_9015}}, trained using the same training data set used for the proposed nonlocal model. As such, this approach yields approximations of $f^*_i(t)$ given by
\begin{equation}\label{eq:NN}
    N\!N(x_i,t;W) \approx f^*_i(t), 
\end{equation}
where $W$ represents the set of weights and biases of the NN, to be learned during training.

\subsection{Problem setting and implementation details}
We consider the two-dimensional problem described in Section \ref{chap:local} with geometry parameters $L=220\sqrt{3}/3$, $l_1=\sqrt{3}/3$, $l_2=1$ and $N=220$, conductivities $\kappa_0 = 1$, $\kappa=0.01$, and hydraulic head $h_0 = 60$. We assume that the particles are injected at 7th unit cell, $C_7$ at time $t=0$ and we track the evolution of their density until time $t=144$.

\begin{figure}[t!]
  \centering
  \includegraphics[scale=0.45]{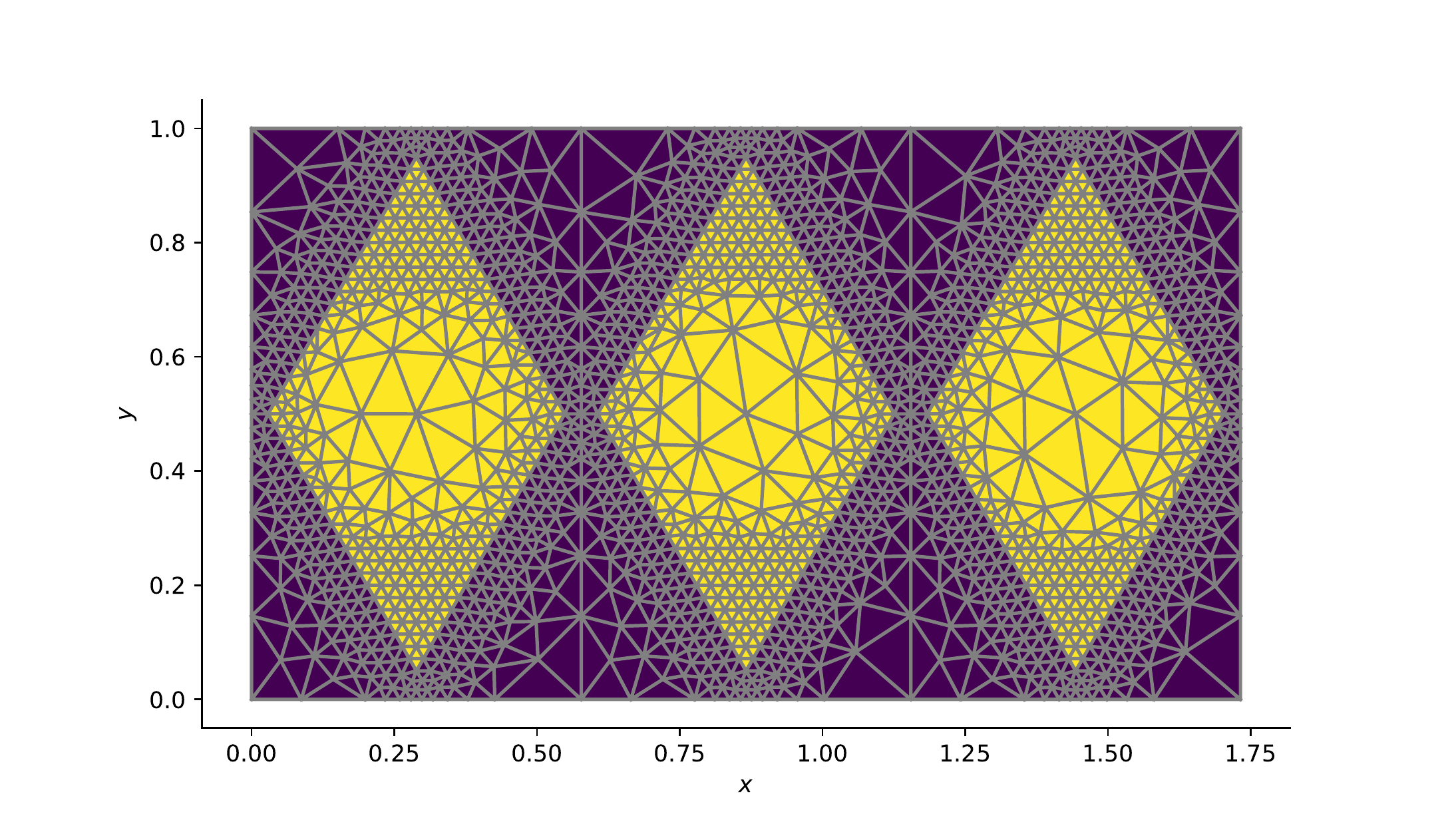}
  \caption{FE mesh of three unit cells.}
  \label{fig:FEMmesh}
\end{figure}
\begin{figure}[t!]
  \centering
  \includegraphics[scale=0.55]{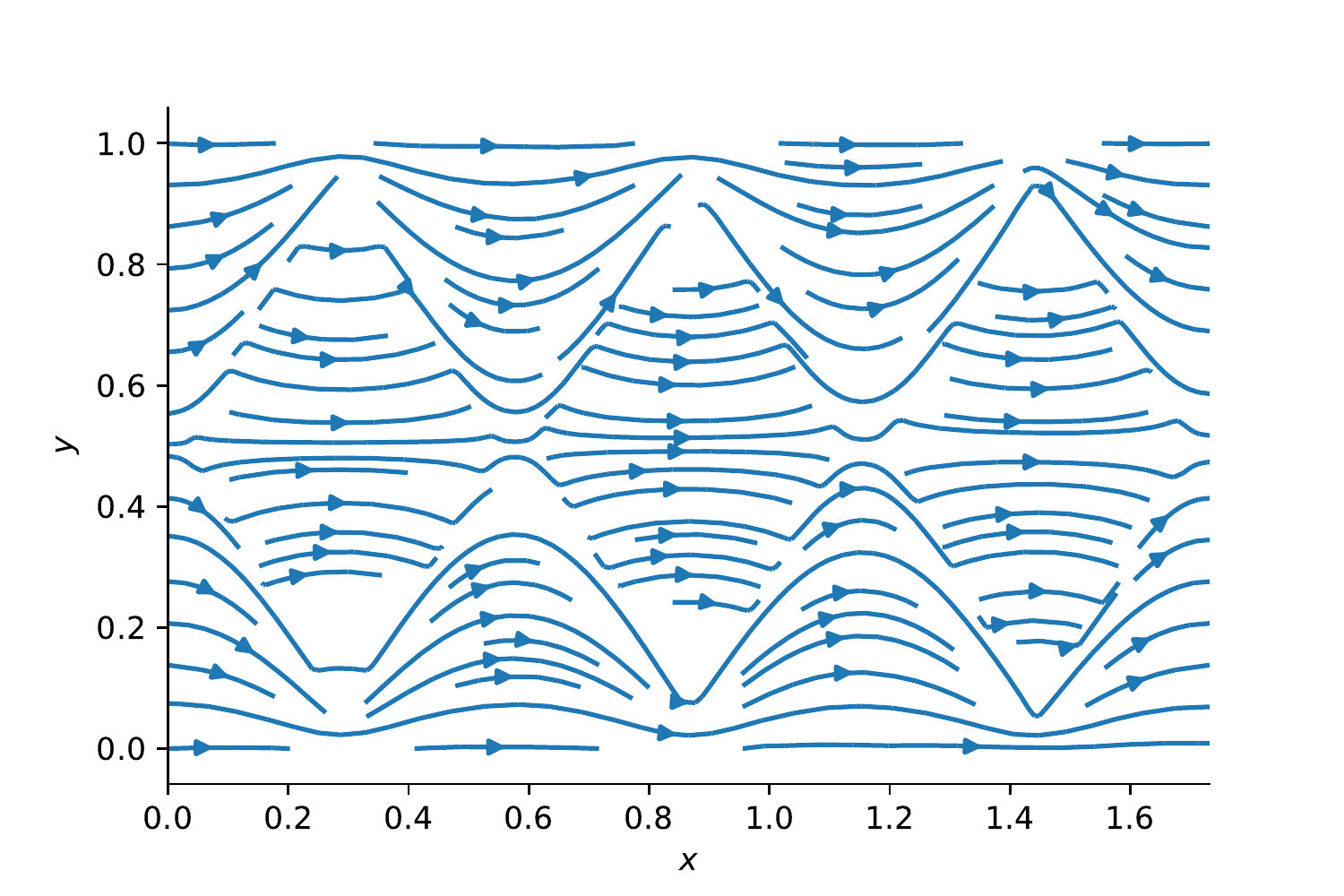}
  \caption{Streamlines in the first 3 unit cells}
  \label{fig:streamline}
\end{figure}
\begin{figure}[t!]
  \centering
  \includegraphics[scale=0.34]{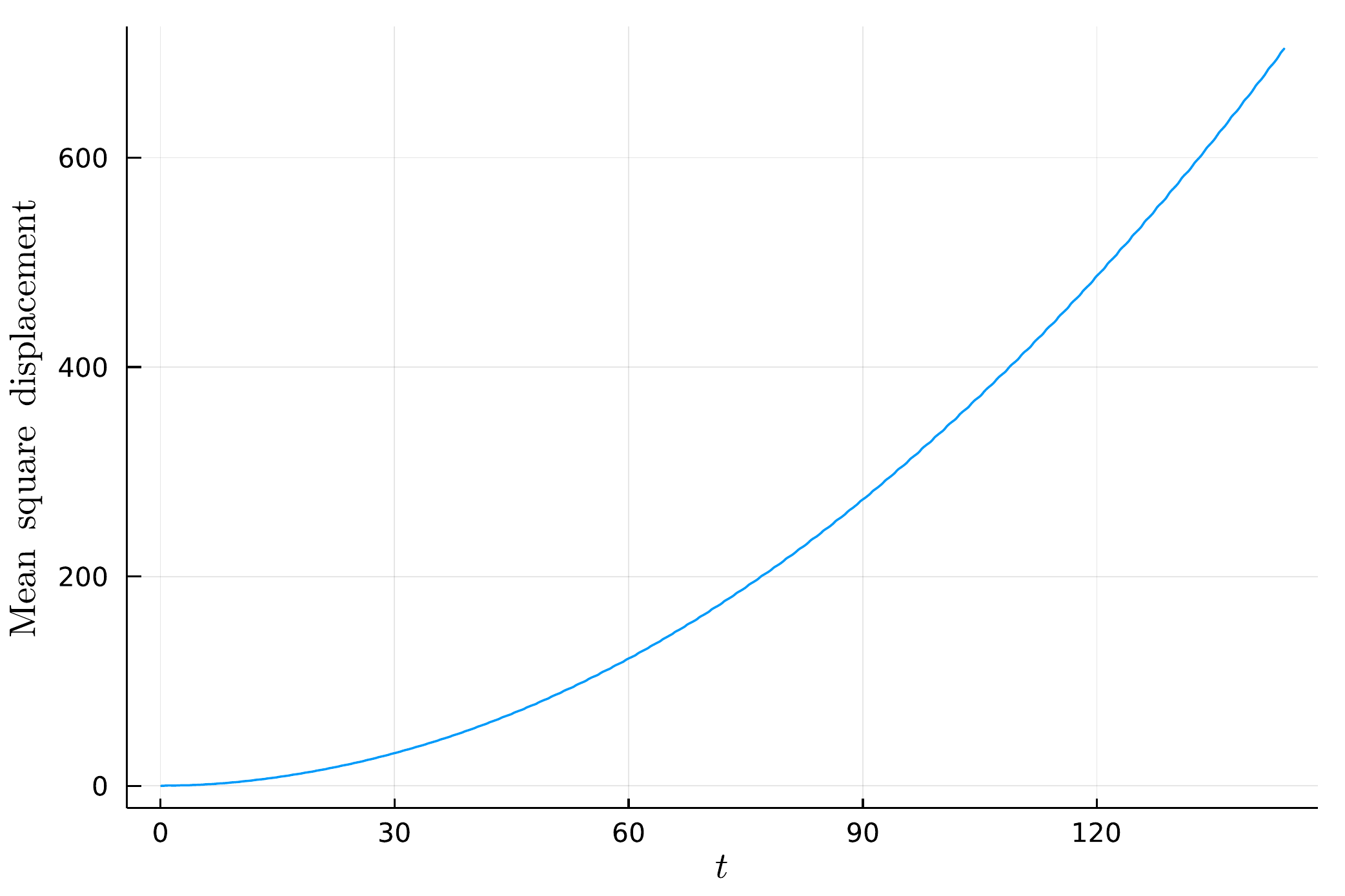}
  \caption{Mean square displacement $\an{x^2} - \an{x}^2$ of the high-fidelity, particle-tracking solution.}
  \label{fig:MSD}
\end{figure}
\begin{figure}[t!]
  \centering
  \includegraphics[scale=0.34]{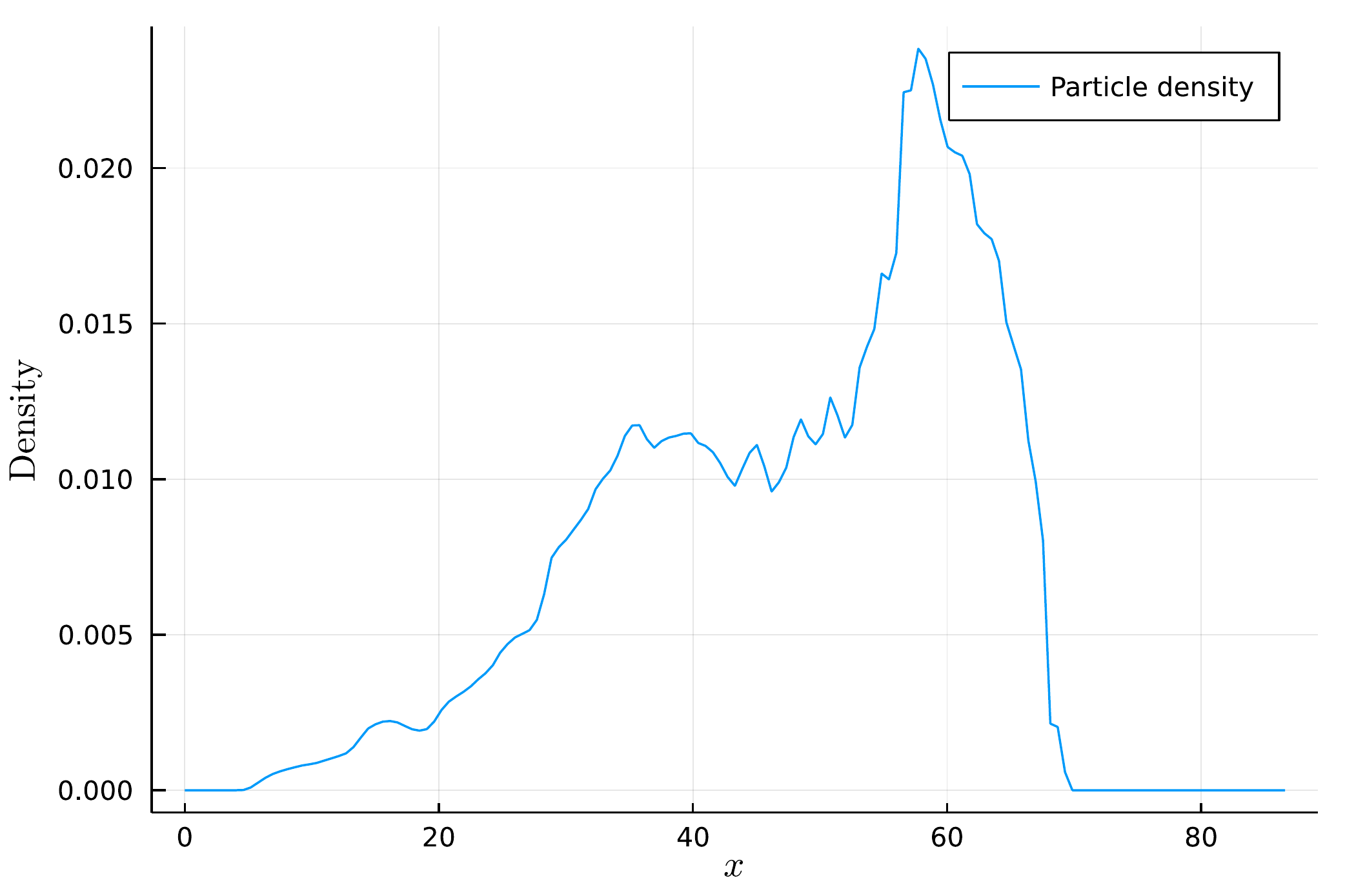}
  \caption{Coarse-grained density profile $\bar c_{H\!F}$ at $t=72$.}
  \label{fig:densitydata}
\end{figure}

\smallskip
For the generation of the high-fidelity data set to be used for training and validation, we solve the flow equation \eref{Poieqn} and the particle transport equation \eref{Lageqn}. Problem \eref{Poieqn} is solved using highly resolved mixed finite element method (FEM). An example of the FEM mesh of three unit cells and the resulting streamlines are reported in Figure~\ref{fig:FEMmesh} and Figure~\ref{fig:streamline} respectively. The particle transport equation \eref{Lageqn} is solved by the ``particle tracking method'' \cite{pollock1988semianalytical} with $22000\times 200$ grid cells, time step $\Delta t=0.01$, and $10^5$ particles. The nonlinear behavior of the corresponding mean squared displacement, reported in Figure~\ref{fig:MSD}, supports our choice of defining the nonlocal diffusion kernel \eref{omegasep} as a function of both space and time. We coarse-grain the high-fidelity density using \eref{upscaling} with $m=20$. An example of the upscaled density at time $t=72$ is reported in Figure~\ref{fig:densitydata}.

For the selection of the training data set we first divide the time interval in two and define the training  interval $[0,T_t]$ and the testing interval $[T_t,144]$ with $0<T_t<144$. Then, we select $n_x=3$ BTCs at locations $\{x_i\}_i=\{ 10,20,30\}$. Thus, the training set is given by
$$
\{ \left(f^*_i(t_j) , x_i \right) \vert \; i=1,\dots,n_x \quad j=1,2,\dots, n_t \},
$$
where, for $T_t = 72$, $n_t=720$. The remaining data set, i.e. the BTCs for $t>T_t$ are used for validation purposes.

\smallskip
The nonlocal diffusion equation \eref{PDeqn} is discretized according to the scheme described in Section \ref{sec:space-discr} with $\delta = N_\delta l_1$, $N_\delta=4$, and $N=220$ (the number of unit cells). We prescribe the initial condition $\bar c(L_1,0) = 1$ by imposing
\begin{align}
    \bar c_i(0) = 1, \quad {\rm for} \; i=7. \notag 
\end{align}

\begin{figure}
  \centering
  \includegraphics[scale=0.4]{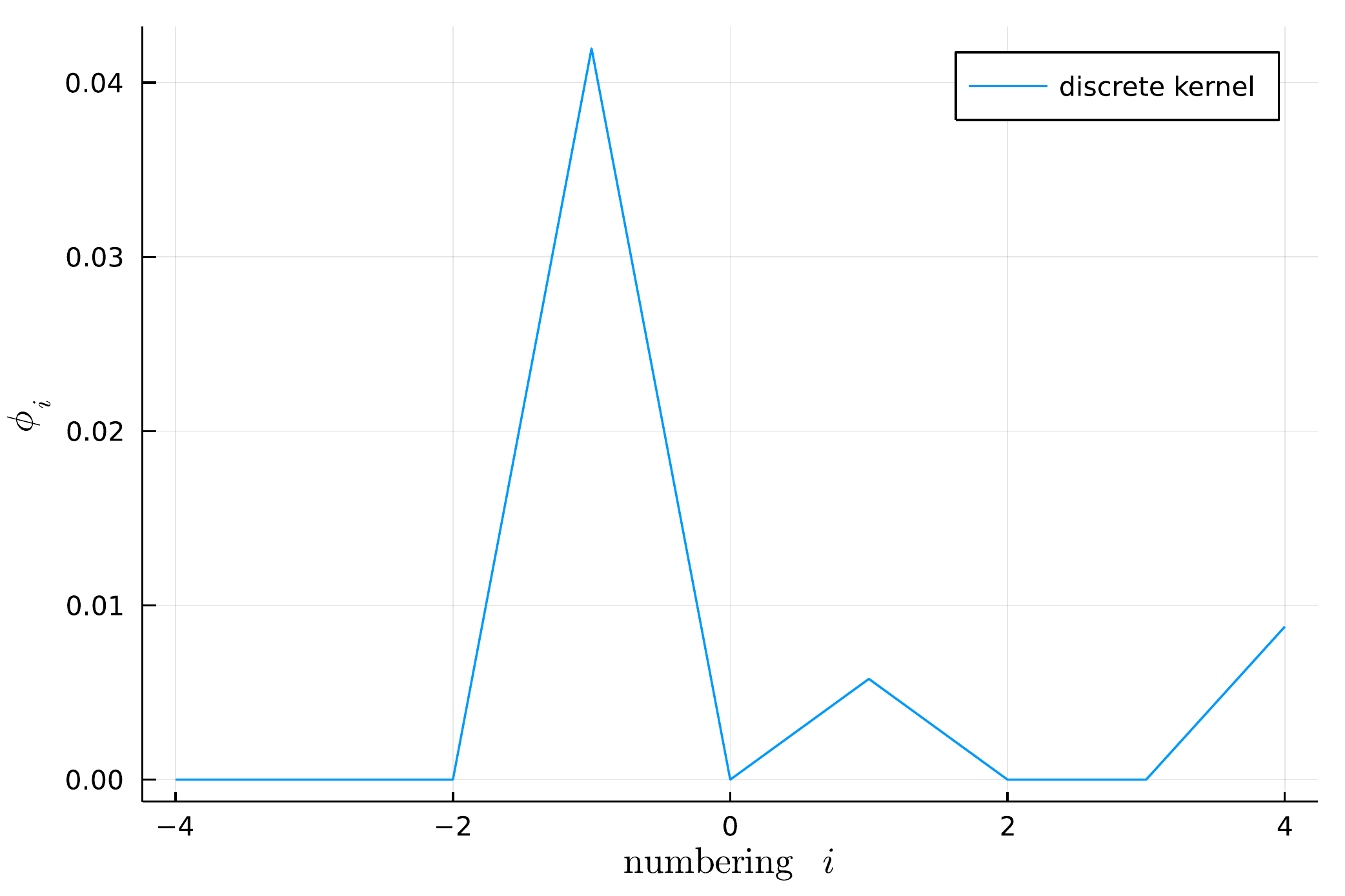}
  \caption{The optimal nonlocal spatial kernel}
  \label{fig:kernel}
\end{figure}

The optimal discrete spatial kernel function $(\nu_\phi,p)$ that minimizes \eqref{eq:discrete-loss} is reported in Figure~\ref{fig:kernel}; here, the optimal power is $p=1.104$.

\subsection{Comparison with alternative PDE models}
\begin{figure}
  \centering
  \includegraphics[scale=0.28]{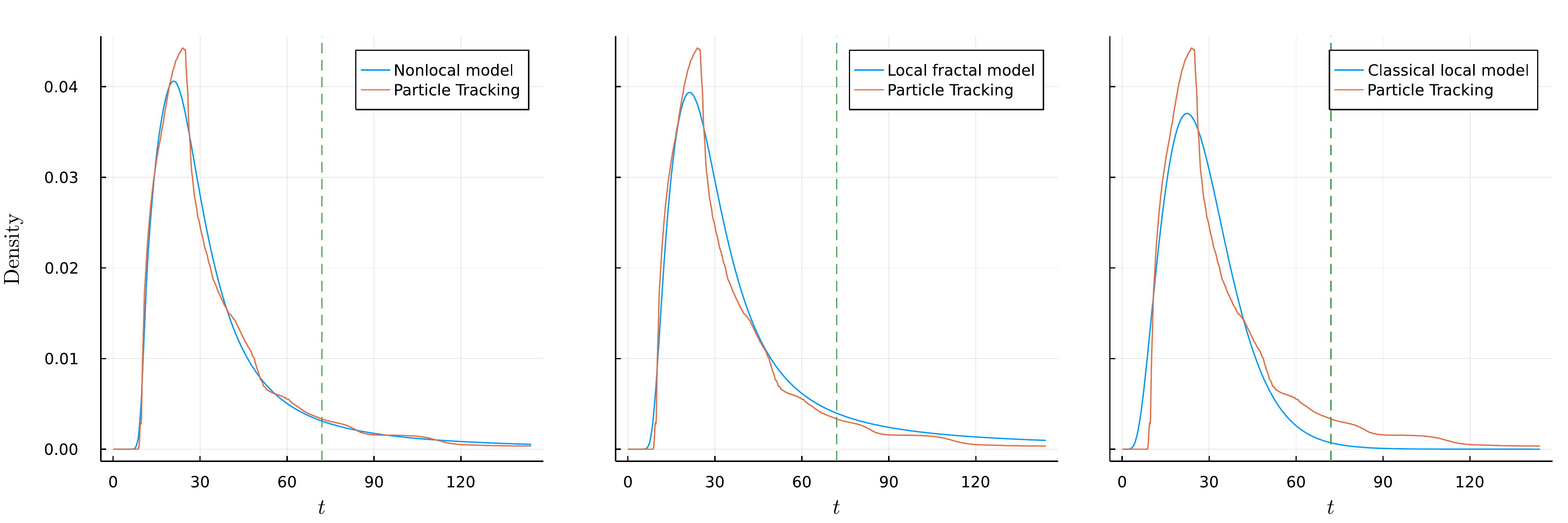}
  \caption{BTC predicted by nonlocal model and local models at location 20 (the green dashed line is used to separate training interval and testing interval)}
  \label{fig:BTC30nnl}
\end{figure}

\begin{figure}
  \centering
  \includegraphics[scale=0.4]{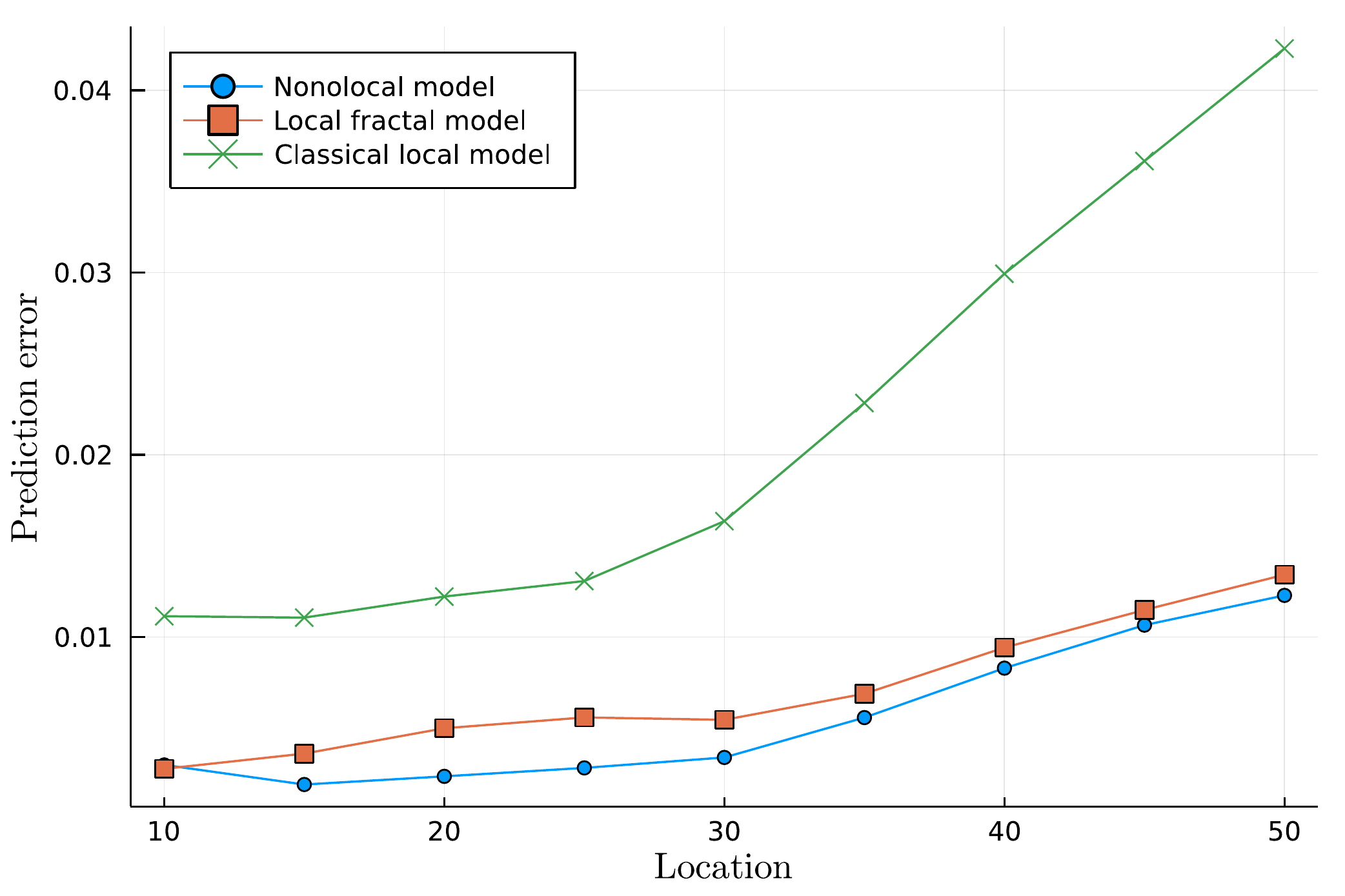}
  \caption{$M\!S\!E$ of predicted BTCs at different locations}
  \label{fig:errornnl}
\end{figure}

In Figure~\ref{fig:BTC30nnl}, we report predictions of the BTC at location $x=20$ obtained by training the nonlocal diffusion model \eqref{eqn:PDdiffeqn}, the fractal model \eqref{eqn:fractaladvdiffeqn} and classical model \eqref{eqn:advdiffeqn}. As a reference, we also report the coarse-grained BTC obtained from the high-fidelity model, $\bar c_{HF}$. We observe that the proposed nonlocal model is in good agreement with the high-fidelity solution, while the predictions obtained with the other PDE models have more pronounced discrepancies when compared to the high-fidelity solution. In particular, we observe that the proposed model has excellent generalization properties; in fact, the predicted BTC, for times that are bigger than the training time $T_t=72$ nearly matches the high-fidelity curve. Instead, it is evident that the classical model fails at generalizing to times $t>T_t$. To further confirm the superiority of the proposed method, we also consider the generalization capability for BTC locations that are different from the ones used for training. We report the $M\!S\!E$ of the predicted BTCs in Figure~\ref{fig:errornnl}: the proposed model always outperforms the alternative models.

\subsection{Comparison with neural network-based surrogates}

According to \eqref{eq:NN}, we consider NN-based surrogates for the coarse-grained BTCs and compare their prediction with the BTCs obtained from the nonlocal solution. We choose a NN with three fully-connected layers, with four nodes per layer, and an output softplus layer to ensure the positivity of the BTC.

\begin{figure}[t!]
  \centering
  \includegraphics[scale=0.4]{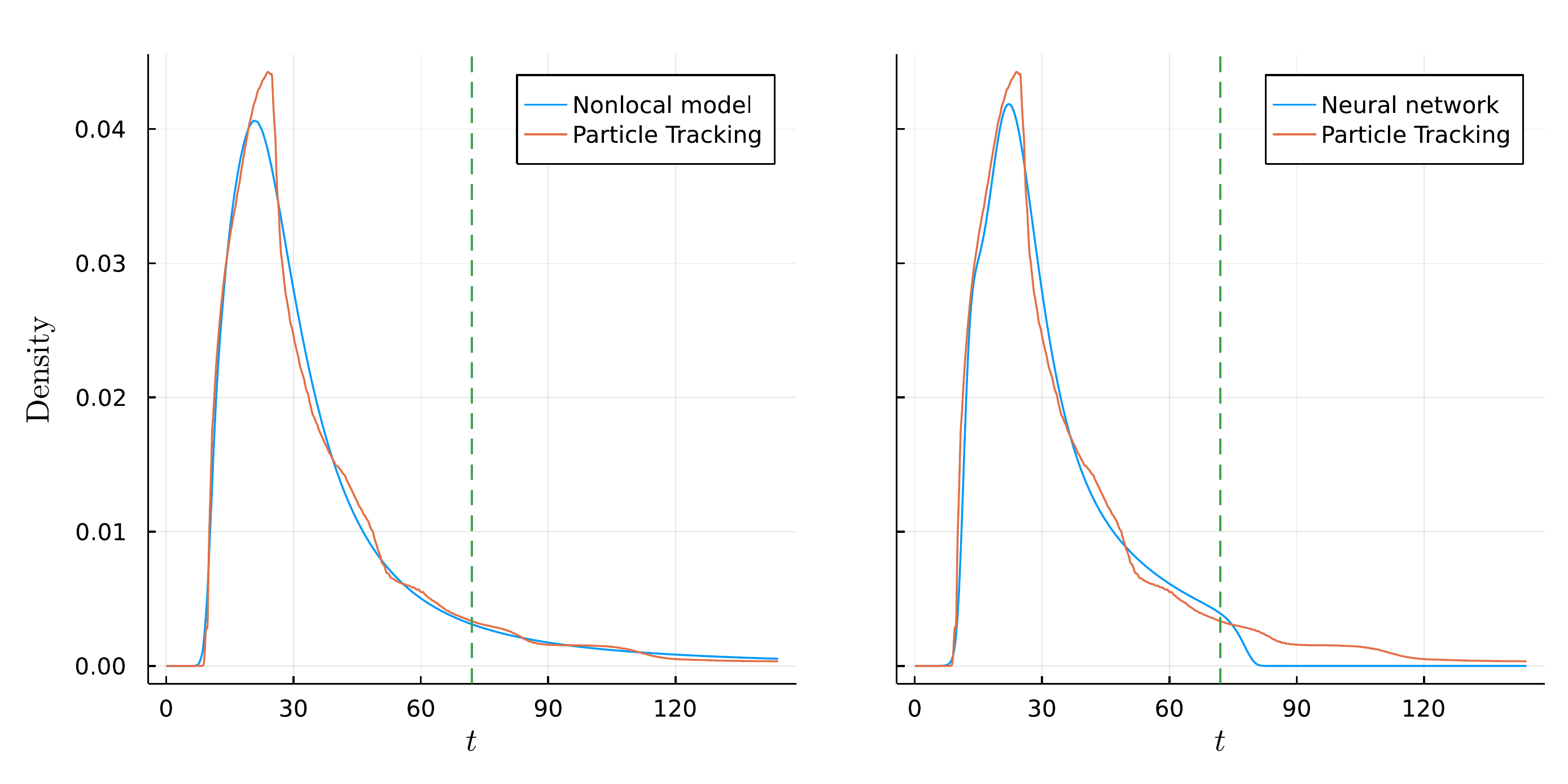}
  \caption{BTCs predicted by nonlocal dynamic kernel and neural network at location 20 (the green dashed line is used to separate training and testing)}
  \label{fig:BTC30nn}
\end{figure}

In the first test case, we train the NN using the same data set described at the beginning of this section. We report the predicted BTC at location $x=20$ together with the corresponding nonlocal prediction in Figure~\ref{fig:BTC30nn}. From these results, as expected, it is evident that the NN has excellent approximation capabilities in the training time interval $[0,72]$, whereas it fails at predicting the BTC at times $t>T_t$. This behavior is due to the fact that the NN over-fitted the training data, compromising its generalization capabilities.

\begin{figure}[t!]
  \centering
  \includegraphics[scale=0.4]{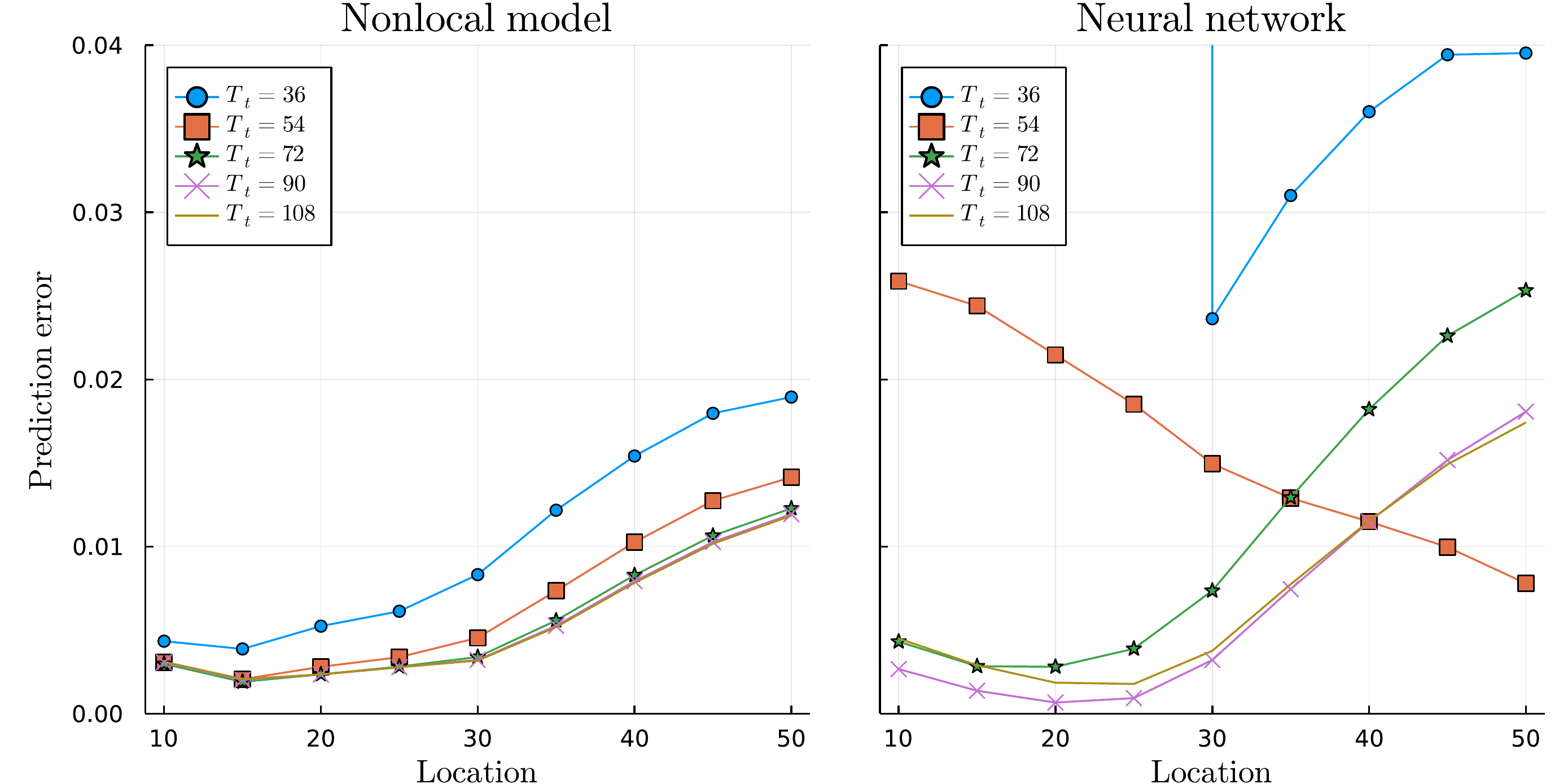}
  \caption{$M\!S\!E$ of predicted BTCs at different location using different size of training data}
  \label{fig:errornn}
\end{figure}

\begin{figure}[t!]
  \centering
  \includegraphics[scale=0.4]{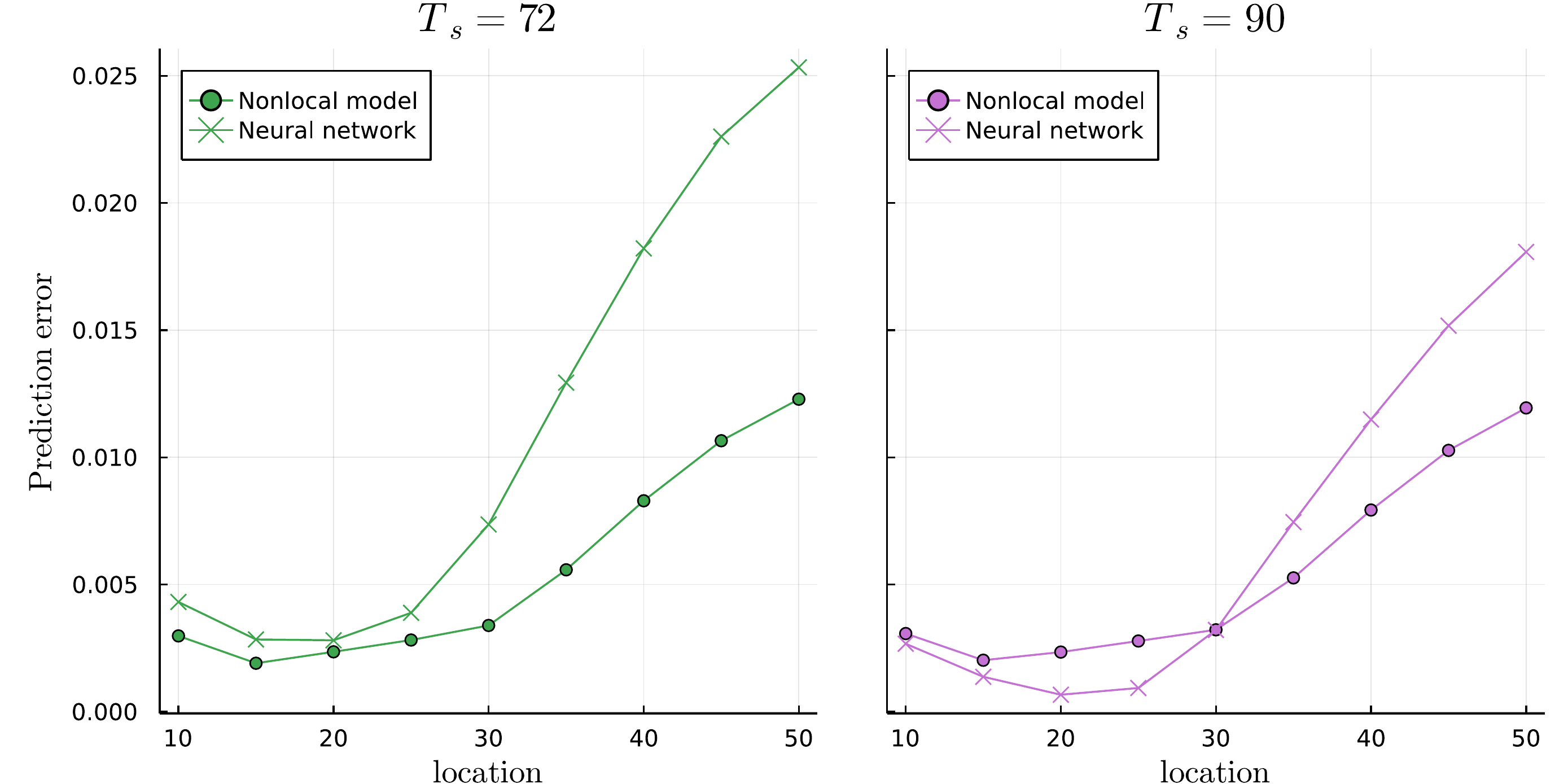}
  \caption{$M\!S\!E$ of predicted BTCs at different location when $T_s = 72$ and $T_s = 90$}
  \label{fig:errornn7290}
\end{figure}

To further compare this approach with the proposed nonlocal model, in a second experiment, we still choose $\{x_i\}_i=\{ 10,20,30 \}$ as training locations, but vary the size of the training interval, i.e. the value of $T_t$. We report the $M\!S\!E$ of the predicted BTCs at several locations (including locations outside the training values) in Figure~\ref{fig:errornn}. We observe that when we decrease the training interval, i.e., decrease the values of $T_t$, the $M\!S\!E$ of the BTCs predicted by the nonlocal model is lower than the $M\!S\!E$ of the NN for locations that are outside the training set. This confirms the good generalization properties of our approach. For locations within the training set the nonlocal and NN prediction errors are similar, but, notably, the superiority of the nonlocal models is more pronounced as $T_t$ decreases. In particular, when $T_s=36$, the neural network fails to make reasonable predictions at locations $\{ 10,15,20,25 \}$. Instead, when $T_s=54$, the NN performs unexpectedly better than the nonlocal model; this is possibly due to the fact that the optimal NN does not over-fit the data.
For the same experiment, with the purpose of providing a visual comparison of the BTCs, we report the $M\!S\!E$ of the predicted BTCs at different locations for $T_t=72$ and $T_t=90$ in Figure~\ref{fig:errornn7290}. These plots confirm our previous statements: when the locations are within the training range, the errors are similar, whereas the nonlocal model outperforms the NN for locations outside the training set. In other words, the NN fails at generalizing beyond the training regime.

\begin{figure}[t!]
  \centering
  \includegraphics[scale=0.4]{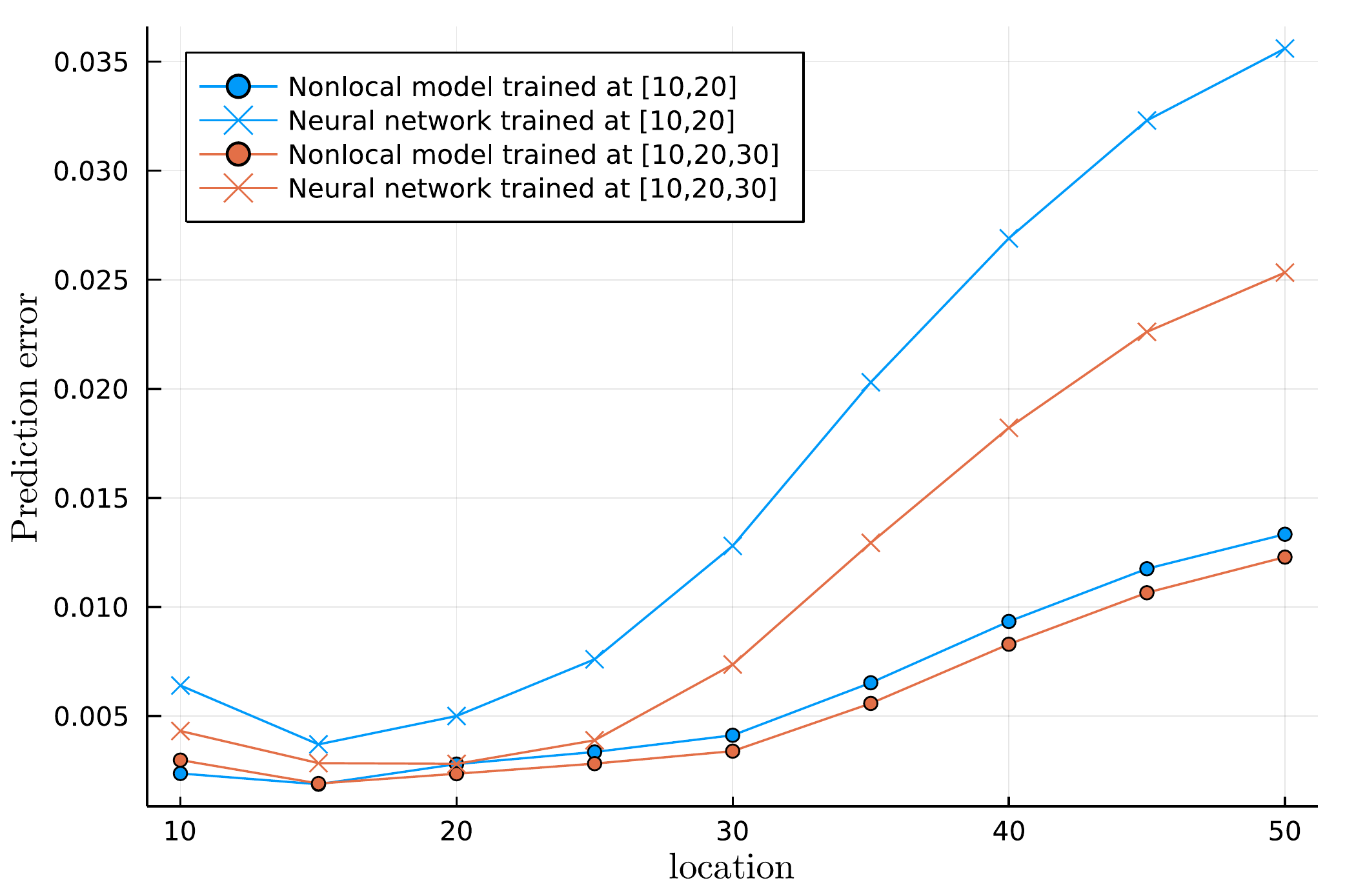}
  \caption{$M\!S\!E$ of predicted BTCs at different location with different choice of training locations}
  \label{fig:error2030}
\end{figure}

In a third experiment, we test the sensitivity of the surrogate to the choice of training locations, while keeping the training interval fixed. We let $T_s=72$ and consider the training locations $\{x_i\}_i=\{10, 20\}$. In Figure~\ref{fig:error2030} we compare the resulting predictions for $M\!S\!E$ with the predictions generated using three training locations and the same training interval. These results indicate, once again, that the proposed approach always outperforms the NN approach, especially for predictions outside the training range.

\section{Conclusion and Future Work}
\label{chap:conclusion}

We presented a data-driven framework where, on the basis of synthetic BTCs, we learn an optimal nonlocal model that describes the evolution of the particle density in a periodic heterogeneous medium. The training data set is generated by means of high-fidelity simulations of small-scale particle transport, that are then upscaled to obtain the coarse-grained BTCs. Specifically, our algorithm learns the nonlocal kernel function such that the corresponding BTCs are as close as possible to the high-fidelity ones. On the basis of numerical comparisons with classical models and black-box NNs, we infer the following.
\begin{enumerate}
    \item The proposed nonlocal dynamic kernel is successful in reproducing anomalous behavior, as it yields a nonlinear mean squared displacement by construction. 
    \item The observed superiority of the proposed model as compared to classical, PDE models confirms that taking into account long-range interactions yields better accuracy for heterogeneous media. Furthermore, both PDE models fail to be predictive outside the training regime (this is particularly pronounced for the classical diffusion model).
    \item Black-box NNs often suffer from overfitting and significantly fail to be predictive outside the training regime. As such, NNs should not be the model of choice in practical situations, where the available data set is sparse in space and limited in time. 
\end{enumerate}

These promising results set the groundwork for more realistic settings and more complex coarse-grained models. In fact, we recall that the high-fidelity model described in Section \ref{chap:local} does not take into account the fine-scale diffusion of the particles which might affect the evolution of the upscaled density and only considers periodic heterogeneity patterns. Also, we only considered spatial nonlocality in the form of a simple nonlocal Laplace operator. Future work includes 1) adding fine-scale diffusion and random heterogeneities to the high-fidelity simulations; 2) the use of field measurements; 3) augmenting the nonlocal model with nonlocality in time; and 4) considering more complex nonlocal operators.


\section*{Acknowledgements}
Xu, D'Elia, Glusa, and Foster are supported by the Sandia National Laboratories (SNL) Laboratory-directed Research and Development program (project 218318). D'Elia is partially supported by the U.S. Department of Energy, Office of Advanced Scientific Computing Research under the Collaboratory on Mathematics and Physics-Informed Learning Machines for Multiscale and Multiphysics Problems (PhILMs) project. SNL is a multimission laboratory managed and operated by National Technology and Engineering Solutions of Sandia, LLC., a wholly owned subsidiary of Honeywell International, Inc., for the U.S. Department of Energy's National Nuclear Security Administration under contract {DE-NA0003525}. This paper, SAND2022-0857 R, describes objective technical results and analysis. Any subjective views or opinions that might be expressed in this paper do not necessarily represent the views of the U.S. Department of Energy or the United States Government.

\clearpage
\appendix
\section{Appendix A}
\label{app:avgv}
Consider a similar Poisson problem as \eref{Poieqn} and \eref{PoieqnBC} but with one single unit cell and unit hydraulic head
\begin{equation}
    \begin{cases}
        \mbf v(x,y) = - \kappa(x,y) \nabla h(x,y) \\
        \nabla \cdot \mbf v(x,y) = 0 
    \end{cases} \quad (x,y) \in [0,l_1]\times[0,l_2] \notag 
\end{equation}
with boundary conditions 
\begin{equation}
    \begin{cases}
        h(0,y) = 1 \\
        h(L,y) = 0 \\
        v_y(x,0) = 0 \\
        v_y(x,l_2) = 0
    \end{cases} \notag 
\end{equation}
Use FEM to solve for the flow field $\mbf v_1(x,y)$. Since particles are injected into the unit cell proportionally to flux in the local problem, the homogenized advection speed $\bar v_1$ can be calculated as the weighted average of advection speed at $x$-direction $v_{1x}(x,y)$ with weights $\vert \mbf v_1(x,y) \vert $ based on classical homogenization techniques \cite{holden1992tensor}:
\begin{align}
    \bar v_1 = \frac{\int_0^{l_1}\int_0^{l_2} \vert \mbf v_1 \vert dy dx}{\int_0^{l_1} \left(\int_0^{l_2} \vert \mbf v_1 \vert dy \right) \frac{\int_0^{l_2} \vert \mbf v_1 \vert dy}{\int_0^{l_2} v_{1x}(x,y) \vert \mbf v_1 \vert dy} dx},
    \label{eqn:avgadv}
\end{align}
To elaborate, the averaging in \eref{avgadv} is done by first arithmetic averaging along $y$-direction with weight $\vert \mbf v_1 \vert$ and then harmonic averaging along $x$-direction. Finally, the effective (homogenized) conductivity $\bar \kappa_x$ can be obtained with
\begin{align}
    \bar \kappa_x = \bar v_1 l_1. \notag 
\end{align}


\bibliographystyle{abbrv}
\bibliography{mybib}
\end{document}